\documentclass[a4paper,12pt]{amsart}
\pdfoutput=1

\usepackage{graphicx}
\usepackage{amsmath}
\usepackage{amssymb}
\usepackage[active]{srcltx}
\usepackage{hyperref}

\newtheorem{thm}{Theorem}%
\theoremstyle{definition}

\theoremstyle{remark}
\theoremstyle{plain}

\def\EE{{\mathbb E}}

\def\HH{{\mathbb H}}

\def\NN{{\mathbb N}}

\def\RR{{\mathbb R}}
\def\TT{{\mathbb T}}

\def\ZZ{{\mathbb Z}}

\def\hatZZ{\widehat{\mathbb Z}}

\def\vecb{{\text{\boldmath$b$}}}

\def\vecm{{\text{\boldmath$m$}}}

\def\vecp{{\text{\boldmath$p$}}}
\def\vecr{{\text{\boldmath$r$}}}

\def\vecx{{\text{\boldmath$x$}}}

\def\vecnull{{\text{\boldmath$0$}}}

\def\scrA{{\mathcal A}}

\def\scrD{{\mathcal D}}

\def\scrF{{\mathcal F}}

\def\scrS{{\mathcal S}}

\def\fC{{\mathfrak C}}

\def\e{\mathrm{e}}
\def\i{\mathrm{i}}

\def\C{\operatorname{C{}}}

\def\SL{\operatorname{SL}}
\def\ASL{\operatorname{ASL}}

\def\meas{\operatorname{meas}}

\def\vol{\operatorname{vol}}

\def\GamG{\Gamma\backslash G}
\def\GamGG{\Gamma_0\backslash G_0}
\def\GamH{\Gamma_H\backslash H}

\def\SLZ{\SL(n+1,\ZZ)}
\def\SLR{\SL(n+1,\RR)}

\def\trans{\,^\mathrm{t}\!}

\title{Fine-scale statistics for the multidimensional Farey sequence}
\author{Jens Marklof}
\address{School of Mathematics, University of Bristol,
Bristol BS8 1TW, U.K.}
\email{j.marklof@bristol.ac.uk}
\date{3 July 2012}

\dedicatory{Dedicated to Friedrich G\"otze on the occasion of his 60th birthday}

\thanks{J.M.\ is supported by a Royal Society Wolfson Research Merit Award, a Leverhulme Trust Research Fellowship and ERC Advanced Grant HFAKT}

\begin{document}

\begin{abstract}
We generalize classical results on the gap distribution (and other fine-scale statistics) for the one-dimensional Farey sequence to arbitrary dimension. This is achieved by exploiting the equidistribution of horospheres in the space of lattices, and the equidistribution of Farey points in a certain subspace of the space of lattices. The argument follows closely the general approach developed by A. Str\"ombergsson and the author [Annals of Math. 172 (2010) 1949--2033].
\end{abstract}

\maketitle

Denote by $\hatZZ^{n+1}$ the set of integer vectors in $\RR^{n+1}$ with relatively prime coefficients, i.e.,
$\hatZZ^{n+1}=\{ \vecm\in\ZZ^{n+1}\setminus\{\vecnull\} : \gcd(\vecm)=1 \}$. The Farey points of level $Q\in\NN$ are defined as the finite set
\begin{equation}\label{fareydef}
	\scrF_Q=\bigg\{ \frac{\vecp}{q} \in[0,1)^n : (\vecp,q)\in\hatZZ^{n+1}, \; 0<q\leq Q \bigg\} .
\end{equation}
The number of Farey points of level $Q$ is asymptotically, for large $Q$,
\begin{equation}\label{asymQ}
	|\scrF_Q| \sim \sigma_Q:=\frac{Q^{n+1}}{(n+1)\,\zeta(n+1)}  .
\end{equation}
In fact, for any bounded set $\scrD\subset[0,1)^n$ with boundary of Lebesgue measure zero and non-empty interior,
\begin{equation}\label{asymQ2}
	|\scrF_Q\cap\scrD| \sim  \vol(\scrD)\; \sigma_Q ,
\end{equation}
which means the Farey sequence is uniformly distributed in $[0,1)^n$.

The objective of the present paper is to understand the fine-scale statistical properties of $\scrF_Q$. To this end, it will be convenient to identify $[0,1)^n$ with the unit torus $\TT^n=\RR^n/\ZZ^n$ via the bijection $[0,1)^n\to\TT^n$, $\vecx\mapsto \vecx+\ZZ^n$. We will consider the following two classical statistical measures of randomness of a deterministic point process: Given $k\in\ZZ_{\geq 0}$ and two test sets $\scrD\subset\TT^n$ and $\scrA\subset\RR^n$, both bounded, with boundary of Lebesgue measure zero and non-empty interior, define
\begin{equation}\label{P}
P_Q(k,\scrD,\scrA)=\frac{\vol\{\vecx\in\scrD : |(\vecx +\sigma_Q^{-1/n} \scrA +\ZZ^n )\cap\scrF_Q| = k\}}{\vol(\scrD)}
\end{equation}
and 
\begin{equation}\label{P0}
P_{0,Q}(k,\scrD,\scrA)=\frac{|\{\vecr\in\scrF_Q\cap\scrD : |(\vecr + \sigma_Q^{-1/n} \scrA +\ZZ^n )\cap\scrF_Q| = k\}|}{|\scrF_Q\cap\scrD|}.
\end{equation}
The scaling of the test set $\scrA$ by a factor $\sigma_Q^{-1/n}$ ensures that the expectation value 
\begin{equation}
\EE P_Q(k,\scrD,\scrA):=\sum_{k=0}^\infty k P_Q(k,\scrD,\scrA)
\end{equation}
is asymptotic to $\vol(\scrA)$ for large $Q$. The quantity $P_{0,Q}(0,\scrD,\scrA)$ is the natural higher dimensional generalization of the gap distribution of sequences in one dimension, which, the case of the Farey sequence for $\scrA=[0,s]$ and $\scrD=\TT$, was calculated by Hall \cite{Hall70}. $P_Q(0,\TT,[0,s])$ corresponds in one dimension to the probability that the distance between a random point on $\TT$ and the nearest element of the sequence is at least $s/2$. An elementary argument shows that in one dimension the density of this distribution is equal to $P_{0,Q}(0,\TT,[0,s])$, see~e.g.~\cite[Theorem 2.2]{nato} and eq.~\eqref{tja} below. The most comprehensive result in one dimension is due to Boca and Zaharescu \cite{Boca05}, who calculate the limiting $n$-point correlation measures. We refer the reader to the survey article \cite{Boca06} for an overview of the relevant literature. 

The results we will discuss here are valid in arbitrary dimension, and will also extend to the distribution in several test sets $\scrA_1,\ldots,\scrA_s$. To keep the notation simple, we will restrict the discussion to one test set; the proofs are otherwise identical, cf.~\cite[Section 6]{partI} for the necessary tools. 

It is evident that the distribution of Farey sequences is intimately linked to the distribution of directions of visible lattice points studied in \cite[Section 2]{partI}. The only difference is in the ordering of the sequence of primitive lattice points and the way they are projected: In the Farey case we take all primitive lattice points in a blow-up of the polytope $\{ (\vecx,y)\in(0,1]^{n+1}: x_j\leq y\}$, draw a line from each lattice point to the origin and and record the intersection of these lines with the hyperplane $\{(\vecx,1):\vecx\in\RR^n\}$. In the case of directions, we take all points in a fixed cone with arbitrary cross-section projected radially onto the unit sphere. Since the cross section of the cone is arbitrary, this yields (by a standard approximation argument) the statistics of primitive lattice points in the blow-up of any star-shaped domain (with boundary of measure zero), which are projected radially onto a suitably chosen hypersurface of codimension one. The proof of a limit distribution for $P_Q(k,\scrD,\scrA)$ for Farey fractions is therefore a corollary of the results of \cite{partI}. 

If the points in $\scrF_Q$ were independent, uniformly distributed random variables in $\TT^n$, we would have, almost surely, convergence to the Poisson distribution:
\begin{equation}
\lim_{Q\to\infty} P_Q(k,\scrD,\scrA) = \lim_{Q\to\infty} P_{0,Q}(k,\scrD,\scrA) = \frac{\vol(\scrA)^k}{k!}\e^{-\vol(\scrA)} \qquad \text{a.s.}
\end{equation}
The $\scrF_Q$ are of course not Poisson distributed. But, as we will see, the limit distributions exist, are independent of $\scrD$, and are given by probability measures on certain spaces of random lattices in $\RR^{n+1}$. The reason for this is as follows. 

Define the matrices
\begin{equation}
h(\vecx)=\begin{pmatrix} 1_n & \trans\vecnull \\ -\vecx & 1 \end{pmatrix},\qquad a(y)=\begin{pmatrix} y^{1/n} 1_n & \trans\vecnull \\ \vecnull & y^{-1} \end{pmatrix}
\end{equation}
and the cone
\begin{equation}
\fC(\scrA)=\{ (\vecx,y)\in\RR^n\times (0,1] : \vecx \in \sigma_1^{-1/n} y\scrA \} \subset \RR^{n+1}.
\end{equation}
Then, for any $(\vecp,q)\in\RR^{n+1}$, 
\begin{equation}
\frac{\vecp}{q}\in \vecx +\sigma_Q^{-1/n} \scrA,\qquad 0<q\leq Q,
\end{equation}
if and only if
\begin{equation}
(\vecp,q) h(\vecx) a(Q) \in \fC(\scrA) .
\end{equation}
Thus, if $Q$ is sufficiently large so that $\sigma_Q^{-1/n}\scrA\subset(0,1]^n$, then
\begin{equation}
|(\vecx + \sigma_Q^{-1/n} \scrA +\ZZ^n )\cap\scrF_Q| =|\hatZZ^{n+1}h(\vecx) a(Q) \cap \fC(\scrA) | .
\end{equation}
This observation reduces the question of the distribution of the Farey sequence to a problem in the geometry of numbers. In particular, \eqref{P} and \eqref{P0} can now be expressed as
\begin{equation}
P_Q(k,\scrD,\scrA)=\frac{\vol\{\vecx\in\scrD : |\hatZZ^{n+1}h(\vecx) a(Q) \cap \fC(\scrA) | = k\}}{\vol(\scrD)}
\end{equation}
and 
\begin{equation}
P_{0,Q}(k,\scrD,\scrA)=\frac{|\{\vecr\in\scrF_Q\cap\scrD : |\hatZZ^{n+1}h(\vecr) a(Q) \cap \fC(\scrA) | = k\}|}{|\scrF_Q\cap\scrD|}.
\end{equation}
Let $G=\SLR$ and $\Gamma=\SLZ$. The quotient $\GamG$ can be identified with the space of lattices in $\RR^{n+1}$ of covolume one. We denote by
$\mu$ the unique right $G$-invariant probability measure on $\GamG$. 
Let furthermore be $\mu_0$ the right $G_0$-invariant probability measure on $\GamGG$, with $G_0=\SL(n,\RR)$ and $\Gamma_0=\SL(n,\ZZ)$.

Define the subgroups
\begin{equation}\label{Hdef}
	H = \bigg\{ M\in G :  (\vecnull,1)M =(\vecnull,1) \bigg\}= \bigg\{ \begin{pmatrix} A  & \trans\vecb \\ \vecnull & 1 \end{pmatrix} : A\in G_0,\; \vecb\in\RR^n \bigg\}
\end{equation}
and
\begin{equation}
	\Gamma_H = \Gamma\cap H  = \bigg\{ \begin{pmatrix} \gamma & \trans\vecm \\ \vecnull & 1 \end{pmatrix} : \gamma\in\Gamma_0,\; \vecm\in\ZZ^n \bigg\} .
\end{equation}
Note that $H$ and $\Gamma_H$ are isomorphic to $\ASL(n,\RR)$ and $\ASL(n,\ZZ)$, respectively. We normalize the Haar measure $\mu_H$ of $H$ so that it becomes a probability measure on $\GamH$. That is,
\begin{equation} \label{siegel2}
d\mu_H(M) = d\mu_0(A)\, d\vecb, \qquad M=\begin{pmatrix} A & \trans\vecb \\ \vecnull & 1 \end{pmatrix}.
\end{equation}

The main ingredient in the proofs of the limit theorems for $P_{0,Q}(k,\scrD,\scrA)$ and $P_Q(k,\scrD,\scrA)$ are the following two equidistribution theorems. The first is the classic equidistribution theorem for closed horospheres of large volume (cf.~\cite[Section 5]{partI} for background and references), the second the equidistribution of Farey points on closed horospheres \cite[Theorem 6]{Marklof10}. In the latter, a key observation is that \cite[eq.~(3.53)]{Marklof10}
\begin{equation}
\Gamma h(\vecr)a(Q) \in \Gamma\backslash\Gamma H a(\tfrac{Q}{q}) \simeq \GamH a(\tfrac{Q}{q}).
\end{equation}

\begin{thm}\label{equiThm0}
For $f:\TT^n\times\GamG\to\RR$ bounded continuous,
\begin{equation}
	\lim_{Q\to\infty} \int_{\TT^n} f\big(\vecx,h(\vecx)a(Q)\big)\, d\vecx  = \int_{\TT^n\times\GamG} f(\vecx,M) \, d\vecx \, d\mu(M) .
\end{equation}
\end{thm}

\begin{thm}\label{equiThm1}
For $f:\TT^n\times\GamG\to\RR$ bounded continuous,
\begin{equation}\label{equiThm1eq}
	\lim_{Q\to\infty} \frac{1}{|\scrF_Q|} \sum_{\vecr\in\scrF_Q} f\big(\vecr,h(\vecr)a(Q)\big)  \\
	= \int_0^1 \int_{\TT^n\times\GamH} f(\vecx, M a(\lambda^{-\frac{1}{n+1}})) \, d\vecx \, d\mu_H(M) \, d\lambda .
\end{equation}
\end{thm}

Both theorems can be derived from the mixing property of the action of the diagonal subgroup $\{a(y)\}_{y\in\RR_{>0}}$. The exponential decay of correlations of this action was exploited by H. Li to calculate explicit rates of convergence \cite{Li11}. One can furthermore generalize Theorem \ref{equiThm1} to general lattices $\Gamma$ in $G$ and non-closed horospheres \cite{HF}.
Theorem \ref{equiThm1} may also be interpreted has an equidistribution theorem for periodic points of the return map of the horocycle flow (in the case $n=1$) to the section 
\begin{equation}\label{trans}
\Gamma\backslash\Gamma H \{ a(y) : y\in\RR_{>1}\}\simeq \GamH  \{ a(y) : y\in\RR_{>1}\}
\end{equation}
which is discussed in \cite{Athreya12}. The identification of \eqref{trans} as an embedded submanifold, which is transversal to closed horospheres of large volume, is central to the proof of Theorem \ref{equiThm1} in \cite{Marklof10}.

By standard probabilistic arguments, the statements of both theorems remain valid if $f$ is a characteristic function of a subset $\scrS\subset\TT^n\times\GamG$ whose boundary has measure zero with respect to the limit measure $d\vecx \, d\mu(M)$ or $d\vecx \, d\mu_H(M) \, d\lambda$, respectively. The relevant set in our application is
\begin{equation}
\scrS = \scrD \times \{ M\in\GamG : |\hatZZ^{n+1}M \cap \fC(\scrA) | \geq k\} .
\end{equation}
The fact that $\scrS$ has indeed boundary of measure zero with respect to $d\vecx \, d\mu(M)$ is proved in \cite[Section 6]{partI}.
We can therefore conclude:
\begin{thm}
Let $k\in\ZZ_{\geq 0}$, and $\scrD\subset\TT^n$, $\scrA\subset\RR^n$ bounded with boundary of Lebesgue measure zero. Then
\begin{equation}
\lim_{Q\to\infty} P_Q(k,\scrD,\scrA)= P(k,\scrA) 
\end{equation}
with
\begin{equation}
P(k,\scrA) = \mu(\{ M\in\GamG : |\hatZZ^{n+1}M \cap \fC(\scrA) | = k\}) ,
\end{equation}
which is independent of the choice of $\scrD$.
\end{thm}

In the second case, we require that the set
\begin{equation}
\{ M\in\GamH : |\hatZZ^{n+1}M \cap \fC_\lambda(\scrA) | \geq k\} , \qquad \fC_\lambda(\scrA):=\fC(\scrA) a(\lambda^{\frac{1}{n+1}}),
\end{equation}
has boundary of measure zero with respect to $\mu_H$, which follows from analogous arguments. With this, we have:

\begin{thm}
Let $k\in\ZZ_{\geq 0}$, and $\scrD\subset\TT^n$, $\scrA\subset\RR^n$ bounded with boundary of Lebesgue measure zero. Then
\begin{equation}\label{Pp}
\lim_{Q\to\infty} P_{0,Q}(k,\scrD,\scrA)= P_0(k,\scrA)= \int_0^1  p_0(k,\fC_\lambda(\scrA))\, d\lambda .
\end{equation}
where
\begin{equation}
p_0(k,\fC) =\mu_H(\{ M\in\GamH : |\hatZZ^{n+1}M \cap \fC) | = k\}),
\end{equation}
which is independent of the choice of $\scrD$.
\end{thm}

In dimension $n\geq 2$, it is difficult to obtain a more explicit description of the limit distributions $P(k,\scrA)$ and $P_0(k,\scrA)$. It is however possible to provide asymptotic estimates for large and small sets $\scrA$ when $k=0$ and $k=1$, see \cite{Athreya09,Strombergsson11}  for general results in this direction. The case of fixed $\scrA$ and large $k$ is discussed in \cite{Marklof00}. 

The geometry of $\GamG$ is significantly simpler in the case $n=1$. This permits the derivation of explicit formulas for the limit distributions in many instances, cf.~\cite{Elkies04,Strombergsson05,partIII}. For example, take $\scrA=(0,s]$, and the cone $\fC_\lambda(\scrA)$ becomes the triangle
\begin{equation}
\Delta_{s,\lambda} = \{ (x_1,x_2)\in\RR^2 : 0<x_1\leq \tfrac{\pi^2}{3} x_2 \lambda s, \; 0<x_2\leq \lambda^{-1/2}  \} ,
\end{equation}
where we have used $\sigma_1=\frac{1}{2\zeta(2)}=\frac{3}{\pi^2}$.
Furthermore $\GamH$ is simply the circle $\TT=\RR/\ZZ$, $\mu_H$ is the standard Lebesgue measure. Hence
\begin{equation}
p_0(k,\Delta_{s,\lambda}) =\meas(\{ x \in \TT :  |\{ (p,q)\in\hatZZ^{2}: (p,px+q)\in\Delta_{s,\lambda}\} | = k\}) .
\end{equation}
It is now a geometric exercise to work out the case $k=0$: With the shorthand $y=\lambda^{1/2}$ and $a=(\tfrac{\pi^2}{3} s)^{-1}$, we deduce
\begin{equation}
p_0(0,\Delta_{s,\lambda}) = 
\begin{cases}
1 & \text{if $y \leq a$}\\
1-\frac{1}{y} +\frac{a}{y^2} & \text{if $a< y \leq a(1-y)^{-1}$}\\
0 & \text{$y > a (1-y)^{-1}$.}\\
\end{cases}
\end{equation}
Solving for $y$, we have in the case $0<a\leq \frac14$
\begin{equation}
p_0(0,\Delta_{s,\lambda}) = 
\begin{cases}
1 & \text{if $y\in[0,a]$}\\
1-\frac{1}{y} + \frac{a}{y^2} & \text{if $y\in[a,\frac12 -\sqrt{\frac14 -a}]\cup [\frac12 + \sqrt{\frac14 -a},1]$ }\\
0 & \text{if $y\in[\frac12 -\sqrt{\frac14 -a},\frac12 + \sqrt{\frac14 -a}]$.}\\
\end{cases}
\end{equation}
For $\frac14<a<1$, we have
\begin{equation}
p_0(0,\Delta_{s,\lambda}) = 
\begin{cases}
1 & \text{if $y\in[0,a]$}\\
1-\frac{1}{y} + \frac{a}{y^2} & \text{if $y\in[a,1],$ }
\end{cases}
\end{equation}
and for $a\geq 1$, we have
\begin{equation}
p_0(0,\Delta_{s,\lambda}) = 1, \qquad y\in[0,1].
\end{equation}
The gap distribution $P_0(0,[0,s])$ is now an elementary integral (recall \eqref{Pp}), which yields
\begin{equation}\label{formula}
P_0(0,[0,s]) = 
\begin{cases}
1 & \text{if $a\in[1,\infty)$}\\
-1+2a-2a\log a & \text{if $a\in[\frac14,1]$}\\
-1+ 2a+2\sqrt{\frac14-a}-4a\log(\frac12+\sqrt{\frac14-a}) & \text{if $a\in[0,\frac14]$.}
\end{cases}
\end{equation}
which reproduces Hall's distribution \cite{Hall70}. The density of this distribution is
\begin{equation}
-\frac{d}{ds} P_0(0,[0,s]) 
=\frac{\pi^2}{3} a^2 \frac{d}{da} P_0(0,[0,s]) 
=
\begin{cases}
0 & \text{if $a\in[1,\infty)$}\\
-\frac{2\pi^2}{3} a^2 \log a & \text{if $a\in[\frac14,1]$}\\
-\frac{4\pi^2}{3} a^2 \log(\frac12+\sqrt{\frac14-a}) & \text{if $a\in[0,\frac14]$,}
\end{cases}
\end{equation}
cf.~\cite[Theorem 2.1]{Boca06}. By \cite[Theorem 2.2]{nato}, we have 
\begin{equation}\label{tja}
-\frac{d}{ds} P(0,[0,s])=P_0(0,[0,s]),
\end{equation}
and hence formula \eqref{formula} yields directly the density of the distribution of the distance to the nearest element. Formula \eqref{formula} was rediscovered in \cite[Lemma 2.6]{Kargaev97}. 

Theorems \ref{equiThm0} and \ref{equiThm1}  reduce  in the case $n=1$ to classic statements in the theory of automorphic forms, with precise bounds on the rate of convergence. Sarnak \cite{Sarnak81} proved Theorem \ref{equiThm0} for test functions $f\in\C_0^\infty$ (infinitely differentiable, compactly supported) that are independent of the first coordinate $x$, and showed that the optimal rate of convergence holds if and only if the Riemann Hypothesis is true (this phenomenon was first pointed out by Zagier \cite{Zagier79}). The reason for the appearance of the Riemann zeros is that the only relevant harmonics in the problem are Eisenstein series $E_{2k}(z,s)$ of even weight $2k$, whose poles are located at the poles of 
\begin{equation}\label{zeta}
\sum_{q=1}^\infty \frac{\varphi(q)}{q^{2s}}=\frac{\zeta(2s-1)}{\zeta(2s)} .
\end{equation}
where $\varphi(s)$ is Euler's totient function and $\zeta(s)$ the Riemann zeta function.

Under the Riemann Hypothesis, Sarnak's rate is significantly better than what one would expected from square-root cancellations---it is the square-root of that. If the test function $f$ depends on $x$ (we assume again $f$ is $\C_0^\infty$), the work of Hejhal \cite{Hejhal00} and Str\"ombergsson \cite{Strombergsson04} shows that the convergence rate slows to the square-root of the horocycle length (or worse) as other terms in the harmonics dominate the error coming from of the Riemann zeros. The object replacing the Eisenstein series in this setting is the Poincar\'e series $P_{m,2k}(z,s)$ of weight $2k$.

The proof Theorem \ref{equiThm1} for $n=1$ on the other hand quickly reduces to estimates of sums over Kloosterman sums. To see this, note first of all that the statement of Theorem  \ref{equiThm1} is equivalent to: For every bounded continuous function $f:\TT\times\Gamma\backslash\HH\to\RR$ (where $\HH$ is the complex upper half plane, on which $\Gamma=\SL(2,\ZZ)$ acts by M\"obius transformations) we have
\begin{equation}\label{equiThm1eq1}
	\lim_{Q\to\infty} \frac{1}{|\scrF_Q|} \sum_{q=1}^Q \sum_{p\in\ZZ_q^\times} f\bigg(\frac{p}{q},\frac{\overline p}{q}+\i \frac{Q^2}{q^2} \bigg)  \\
	= \int_0^\infty\int_0^1\int_0^1  f(x, u+\i v) \, dx \, du \, \frac{dv}{v^2} .
\end{equation}
Here $\ZZ_q^\times$ denotes the multiplicative group of invertible residues mod $q$, and $\overline p$ is the inverse of $p\bmod q$. One way of proving \eqref{equiThm1eq1} is to expand $f\in\C_0^\infty$ in its harmonics (Fourier series in $x$ and $y$ and Mellin transform in $\lambda$) which leads to Selberg's Kloosterman zeta function 
\begin{equation}
Z_{m_1,m_2}(s) = \sum_{q=1}^\infty \frac{K(m_1,m_2,q)}{q^{2s}}
\end{equation}
with the Kloosterman sum 
\begin{equation}
K(m_1,m_2,q) = \sum_{p\in\ZZ_q^\times} \e^{2\pi\i (m_1p+m_2\overline p)/q} . 
\end{equation}
As in the case of the equidistribution of closed horocycles, where the asymptotics was determined by the poles of the Eisenstein and Poincar\'e series, the poles of $Z_{m_1,m_2}(s)$ now determine the asymptotics of \eqref{equiThm1eq1}. Note that $Z_{0,m_2}(s)$ are precisely the Fourier coefficients of $E_0(z,s)$ and, as already understood by Selberg, $Z_{m_1,m_2}(s)$ is the $m_2$-th Fourier coefficient of the Poincar\'e series $P_{m_1,0}(z,s)$. Hence the appearance of the Riemann hypothesis in the error term of Theorem \ref{equiThm1} mirrors exactly the situation in Theorem \ref{equiThm0}.

\end{document}